\documentclass{article}
\usepackage{amsmath}
\usepackage{amssymb}
\usepackage{verbatim}

\newtheorem{theorem}{Theorem}
\newtheorem{acknowledgement}[theorem]{Acknowledgement}

\newtheorem{definition}[theorem]{Definition}

\newtheorem{proposition}[theorem]{Proposition}
\newtheorem{remark}[theorem]{Remark}

\input{tcilatex}
\begin{document}

\title{Towards a Morse theory on Banach spaces via ultrafunctions }
\author{Vieri Benci\thanks{
Dipartimento di Matematica, Universit\`{a} degli Studi di Pisa, Via F.
Buonarroti 1/c, 56127 Pisa, ITALY; e-mail: \texttt{benci@dma.unipi.it }and%
\texttt{\ }Centro Linceo Interdisciplinare "Beniamino Segre"},\thanks{%
Visiting UFRJ} \and Isaia Nisoli\thanks{%
Instituto de Matem\'{a}tica - UFRJ Av. Athos da Silveira Ramos 149, Centro
de Tecnologia - Bloco C Cidade Universit\~{A}\textexclamdown ria - Ilha do
Fund\~{a}o. Caixa Postal 68530 21941-909 Rio de Janeiro - RJ - Brasil
Email:nisoli@im.ufrj.br}}
\maketitle

\begin{abstract}
Morse Theory on Banach spaces would be a useful tool in nonlinear analysis
but its development is hindered by many technical problems. In this paper we
present an approach based on a new notion of generalized functions called
\textquotedblleft ultrafunctions\textquotedblright\ which solves some of the
technical questions involved.

\medskip \medskip


\medskip

\end{abstract}

\tableofcontents

\section{Introduction}

In this paper we start a study of Morse Theory on Banach spaces using the
theory of Ultrafunctions \cite{belu2013,milano,beyond,topology}; the \textbf{%
ultrafunctions} are a new notion of generalized functions based on the
general ideas of Non Archimedean Mathematics (NAM) of Non Standard Analysis (NSA).

Based on our experience NAM allows to construct models of the physical world
in a more elegant and simpler way, in many circumstances. Contrary to the
common belief, the ideas behind NSA and NMA date backs to the years around
1870's, when it was investigated by mathematicians such as Du Bois-Reymond,
Veronese, Hilbert and Levi-Civita. Since then its development stopped, until
the '60s when Abraham Robinson presented his Non Standard Analysis (NSA).
For a historical analysis of these facts we refer to Ehrlich \cite{el06} and
to Keisler \cite{keisler76} for a very clear exposition of NSA.

Ultrafunctions are a particular class of functions based on a superreal
field $\mathbb{R}^{\ast }\supset \mathbb{R}$. More exactly, to any
continuous function $f:\mathbb{R}^{N}\rightarrow \mathbb{R}$, we associate
in a canonical way an ultrafunction $\widetilde{f}:\left( \mathbb{R}^{\ast
}\right) ^{N}\rightarrow \mathbb{R}^{\ast }$ which extends $f$; the
ultrafunctions are many more than the functions and among them we can find
solutions of functional equations which do not have any solutions among the
real functions or the distributions.

Also, the theory of ultrafuctions allows to overcome some difficulties
of Morse Theory in Banach spaces.

Many authors have been working on the adaptation of Morse Theory on Banach
spaces \cite{Chang83,Chang93,Chang98,MerPalm,Uhl72}, but many problems
arise: a really important one is the difficulty in defining what a (weakly)
nondegenerate critical point is and how to define its Morse index, since any
critical point of a $C^{2}$ functional on a Banach space is degenerate and
it is not possible to apply the generalized Morse Lemma (for a reference on
the generalized Morse Lemma see \cite{GroMe}).

In recent times, a lot of delicate work has been done in this direction,
developing extremely refined tools and techniques to study problems in
nonlinear analysis \cite%
{CinCarMartVan,CD05,CinDe09,CinDeSciCGpLap,CLV05,CV03,CV06,CV07,CV09,CV09b,CinVanVis11,CinVanVisMul,Lan}%
. Our approach is totally different, we avoid many of the difficulties
involved in the definitions by using the properties of hyperfinite function
spaces.

We believe that the flexibility of the ultrafunction approach can be
fruitful for the development of the Theory. In this paper we present a
foundational basis for this theory; other articles dealing with applications
are to follow.

\bigskip 

\begin{acknowledgement}
The first author wishes to thank the Federal University of Rio de Janeiro
(UFRJ) for the invitation and hospitality.
\end{acknowledgement}

\subsection{Notation}

We fix some notation. Since this paper does not deal with application, 
we use some function spaces as model spaces for the theory; 
let $\Omega $\ be a subset of $\mathbb{R}^{N}$:

\begin{itemize}
\item $\mathcal{C}\left( \Omega \right) $ denotes the set of real continuous
functions defined on $\Omega$;

\item $\mathcal{C}_{0}\left( \overline{\Omega }\right) $ denotes the set of
real continuous functions on $\overline{\Omega }$ which vanish on $\partial
\Omega$;

\item $\mathcal{C}^{k}\left( \Omega \right) $ denotes the set of functions
defined on $\Omega \subset \mathbb{R}^{N}$ which have continuous derivatives
up to the order $k$;

\item $\mathcal{C}_{0}^{k}\left( \overline{\Omega }\right) =\mathcal{C}%
^{k}\left( \overline{\Omega }\right) \cap \mathcal{C}_{0}\left( \overline{%
\Omega }\right) $;

\item $\mathcal{D}\left( \Omega \right) $ denotes the set of the infinitely
differentiable functions with compact support defined on $\Omega \subset 
\mathbb{R}^{N}$;

\item $L^2\left( \Omega \right) $ denotes the set of square integrable functions
on $\Omega$.
\end{itemize}

\section{Preliminary notions}

In this section we present some background material necessary to follow the
following part. We underline that this material is not original but we cite
it in order to make the article (almost) self contained. We refer to \cite%
{belu2013,milano,beyond,topology} for a more detailed treatment.

\subsection{Non Archimedean Fields\label{naf}}

Here, we recall the basic definitions and facts regarding non-Archimedean
fields. In the following, ${\mathbb{K}}$ will denote an ordered field. We
recall that such a field contains (a copy of) the rational numbers. Its
elements will be called numbers.

\begin{definition}
Let $\mathbb{K}$ be an ordered field. Let $\xi \in \mathbb{K}$. We say that:

\begin{itemize}
\item $\xi $ is infinitesimal if, for all positive $n\in \mathbb{N}$, $|\xi
|<\frac{1}{n}$;

\item $\xi $ is finite if there exists $n\in \mathbb{N}$ such as $|\xi |<n$;

\item $\xi $ is infinite if, for all $n\in \mathbb{N}$, $|\xi |>n$
(equivalently, if $\xi $ is not finite).
\end{itemize}
\end{definition}

\begin{definition}
An ordered field $\mathbb{K}$ is called Non-Archimedean if it contains an
infinitesimal $\xi \neq 0$.
\end{definition}

It's easily seen that all infinitesimal are finite, that the inverse of an
infinite number is a nonzero infinitesimal number, and that the inverse of a
nonzero infinitesimal number is infinite.

\begin{definition}
A superreal field is an ordered field $\mathbb{K}$ that properly extends $%
\mathbb{R}$.
\end{definition}

It is easy to show, due to the completeness of $\mathbb{R}$, that there are
nonzero infinitesimal numbers and infinite numbers in any superreal field.
Infinitesimal numbers can be used to formalize a new notion of ``closeness'':

\begin{definition}
\label{def infinite closeness} We say that two numbers $\xi, \zeta \in {%
\mathbb{K}}$ are infinitely close if $\xi -\zeta $ is infinitesimal. In this
case, we write $\xi \sim \zeta $.
\end{definition}

Clearly, the relation ``$\sim $'' of infinite closeness is an equivalence
relation.

\begin{theorem}
If $\mathbb{K}$ is a superreal field, every finite number $\xi \in \mathbb{K}
$ is infinitely close to a unique real number $r\sim \xi $, called the 
\textbf{shadow} or the \textbf{standard part} of $\xi $.
\end{theorem}

Given a finite number $\xi $, we denote its shadow as $sh(\xi )$, and we put 
$sh(\xi )=+\infty $ ($sh(\xi )=-\infty $) if $\xi \in \mathbb{K}$ is a
positive (negative) infinite number.\newline

\begin{definition}
Let $\mathbb{K}$ be a superreal field, and $\xi \in \mathbb{K}$ a number.
The \label{def monad} monad of $\xi $ is the set of all numbers that are
infinitely close to it:%
\begin{equation*}
\mathfrak{m}\mathfrak{o}\mathfrak{n}(\xi )=\{\zeta \in \mathbb{K}:\xi \sim
\zeta \},
\end{equation*}%
and the galaxy of $\xi $ is the set of all numbers that are finitely close
to it: 
\begin{equation*}
\mathfrak{gal}(\xi )=\{\zeta \in \mathbb{K}:\xi -\zeta \ \text{is\ finite}\}
\end{equation*}
\end{definition}

By definition, it follows that the set of infinitesimal numbers is $%
\mathfrak{mon}(0)$ and that the set of finite numbers is $\mathfrak{gal}(0)$.

\subsection{The $\Lambda $-limit\label{OL}}

In this section we will introduce a particular superreal field $\mathbb{K}$
and we will analyze its main properties by means of $\Lambda $-theory, in
particular by means of the notion of $\Lambda $-limit (for complete proofs
and for further properties of the $\Lambda $-limit, the reader is referred
to \cite{ultra,belu2013,milano,beyond,topology}).

We recall that the superstructure on $\mathbb{R}$ is defined as follows:

\begin{equation*}
\mathbb{U}=\dbigcup_{n=0}^{\infty }\mathbb{U}_{n}
\end{equation*}%
where $\mathbb{U}_{n}$ is defined by induction as follows:%
\begin{eqnarray*}
\mathbb{U}_{0} &=&\mathbb{R}\text{;} \\
\mathbb{U}_{n+1} &=&\mathbb{U}_{n}\cup \mathcal{P}\left( \mathbb{U}%
_{n}\right) .
\end{eqnarray*}%
Here $\mathcal{P}\left( E\right) $ denotes the power set of $E.$ Identifying
the couples with the Kuratowski pairs and the functions and the relations
with their graphs, it follows that{\ }$\mathbb{U}$ contains almost every
usual mathematical object. Now, we set 
\begin{equation*}
\mathfrak{L}=\mathcal{P}_{\omega }(\mathbb{U}),
\end{equation*}
and we will refer to $\mathfrak{L}$ as the \textquotedblleft parameter
space\textquotedblright . Clearly $\left( \mathfrak{L},\subset \right) $ is
a directed set\footnote{%
We recall that a directed set is a partially ordered set $(D,\prec )$ such
that, $\forall a,b\in D,\ \exists c\in D$ such that 
\begin{equation*}
a\prec c\ \ \text{and}\ \ b\prec c.
\end{equation*}%
}. We add at $\mathfrak{L}$ one point at infinity $\Lambda $ and we define
the following family of neighborhoods of infinity:%
\begin{equation*}
\left\{ \Lambda \cup Q\ |\ Q\in \mathcal{U}\right\}
\end{equation*}%
where $\mathcal{U}$ is a fine ultrafilter on $\mathfrak{L,}$ namely it is a
filter such that

\begin{itemize}
\item if $A\cup B=\mathfrak{L}$, then 
\begin{equation}
A\in \mathfrak{\mathcal{U}\ }\text{or }B\in \mathcal{U};  \label{quaqua}
\end{equation}

\item $\forall \lambda _{0}\in \mathfrak{L}$, $\left\{ \lambda \in \mathfrak{%
L}\ |\ \lambda _{0}\subset \lambda \right\} \in \mathcal{U}$
\end{itemize}

A function $\varphi :D\rightarrow E$ defined on a directed set will be
called \textit{net }(with values in $E$). If $\varphi _{\lambda }$ is a real
net, we have that 
\begin{equation*}
\underset{\lambda \rightarrow \Lambda }{\lim }\varphi _{\lambda }=L
\end{equation*}%
if and only if 
\begin{equation}
\forall \varepsilon >0,\text{ }\exists Q\in \mathcal{U}\text{\ such that, }%
\forall \lambda \in Q,\ \left\vert \varphi _{\lambda }-L\right\vert
<\varepsilon .  \label{limite}
\end{equation}%
We will refer to the sets in $Q$ as \textbf{qualified sets}.

Notice that this topology on $\mathfrak{L}\cup \left\{ \Lambda \right\} $
satisfies this interesting property:

\begin{proposition}
\label{nino}If the net $\varphi _{\lambda }$ has a converging subnet, then
it is a \textbf{converging} net.
\end{proposition}

\textbf{Proof}: Suppose that the net $\varphi _{\lambda }$ has a converging
subnet to $L\in \mathbb{R}$. We fix $\varepsilon >0$ arbitrarily and we have
to prove that $Q_{\varepsilon }\in \mathcal{U}$ where%
\begin{equation*}
Q_{\varepsilon }=\left\{ \lambda \in \mathfrak{L}\ |\ \left\vert \varphi
_{\lambda }-L\right\vert <\varepsilon \right\} .
\end{equation*}%
We argue indirectly and we assume that 
\begin{equation*}
Q_{\varepsilon }\notin \mathcal{U}
\end{equation*}%
Then, by (\ref{quaqua}), $N=\mathfrak{L}\backslash \left( Q_{\varepsilon
}\cap E\right) \in \mathcal{U}$ and hence%
\begin{equation*}
\forall \lambda \in N,\ \left\vert \varphi _{\lambda }-L\right\vert \geq
\varepsilon ,
\end{equation*}
This contradict the fact that $\varphi _{\lambda }$ has a subnet which
converges to $L.$

$\square $

\bigskip

We have the following result:

\begin{theorem}
\label{nuovo}\textit{There exists a superreal field} $\mathbb{K}\supset 
\mathbb{R}$\textit{\ a Hausdorff topology on the space }$\left( \mathfrak{L}%
\times \mathbb{R}\right) \cup \mathbb{K}$ \textit{such that }

\begin{enumerate}
\item \textit{Every net }$\varphi :\mathfrak{L}\times \mathbb{R}\rightarrow 
\mathbb{R}$\textit{\ has a unique limit } 
\begin{equation*}
L=\lim_{\lambda \rightarrow \Lambda }\left( \lambda ,\varphi (\lambda
)\right) .
\end{equation*}%
\textit{Moreover we assume that every}\emph{\ }$\xi \in \mathbb{K}$\textit{\
is the limit\ of some net }$\varphi :\mathfrak{L}\times \mathbb{R}%
\rightarrow \mathbb{R}$\emph{.}

\item If $r\in \mathbb{R}$ 
\begin{equation*}
\lim_{\lambda \rightarrow \Lambda }\left( \lambda ,r\right) =r.
\end{equation*}

\item \emph{\ }\textit{For all }$\varphi ,\psi :\mathfrak{L}\rightarrow 
\mathbb{R}$\emph{: }%
\begin{eqnarray*}
\lim_{\lambda \rightarrow \Lambda }\left( \lambda ,\varphi (\lambda )\right)
+\lim_{\lambda \rightarrow \Lambda }\left( \lambda ,\psi (\lambda )\right)
&=&\lim_{\lambda \uparrow \Lambda }\left( \lambda ,\varphi (\lambda )+\psi
(\lambda )\right) ; \\
\lim_{\lambda \rightarrow \Lambda }\left( \lambda ,\varphi (\lambda )\right)
\cdot \lim_{\lambda \rightarrow \Lambda }\left( \lambda ,\psi (\lambda
)\right) &=&\lim_{\lambda \rightarrow \Lambda }\left( \lambda ,\varphi
(\lambda )\cdot \psi (\lambda )\right) .
\end{eqnarray*}
\end{enumerate}
\end{theorem}

\textbf{Idea of the proof:} The proof of this theorem is in \cite{topology}.
We now will sketch it for the sake of the reader. We set%
\begin{equation*}
I=\left\{ \varphi \in \mathfrak{F}\left( \mathfrak{L},\mathbb{R}\right) \ |\
\varphi (x)=0\ \text{in a qualified set}\right\} .
\end{equation*}%
It is not difficult to prove that $I$ is a maximal ideal in $\mathfrak{F}%
\left( \mathfrak{L},\mathbb{R}\right) ;$ then%
\begin{equation*}
\mathbb{K}:=\frac{\mathfrak{F}\left( \mathfrak{L},\mathbb{R}\right) }{I}
\end{equation*}%
is a field. In the following, we shall identify a real number $c\in \mathbb{R%
}$ with the equivalence class of the constant net $\left[ c\right] _{I}.$

Now, we equip $\left( \mathfrak{L}\times \mathbb{R}\right) \cup \mathbb{K}$
with the following topology $\tau $. A basis of neighborhoods of $\left[
\varphi \right] _{I}$ is given by 
\begin{equation*}
N_{\varphi ,Q}:=\left\{ \left( \lambda ,\varphi (\lambda )\right) \mid
\lambda \in Q\right\} \cup \left\{ \left[ \varphi \right] _{I}\right\} ,\ \
Q\in \mathcal{U}.
\end{equation*}

$\square $

From now on, in order to simplify the notation we will write%
\begin{equation*}
\lim_{\lambda \uparrow \Lambda }\varphi (\lambda ):=\lim_{\lambda
\rightarrow \Lambda }\left( \lambda ,\varphi (\lambda )\right) ,
\end{equation*}%
and we call it $\Lambda $-limit.

\subsection{Natural extension of sets and functions}

The notion of $\Lambda $-limit can be extended to sets and functions in the
following way:

\begin{definition}
\label{limito}Let $E_{\lambda },$ $\lambda \in \mathfrak{L},$ be a family of
sets in $\mathbb{R}^{N}.$ We pose%
\begin{equation*}
\lim_{\lambda \uparrow \Lambda }\ E_{\lambda }:=\left\{ \lim_{\lambda
\uparrow \Lambda }\psi (\lambda )\ |\ \psi (\lambda )\in E_{\lambda
}\right\} .
\end{equation*}%
A set which is a $\Lambda $-\textit{limit\ is called \textbf{internal}.} In
particular if, $\forall \lambda \in \mathfrak{L,}$ $E_{\lambda }=E,$ we set $%
\lim_{\lambda \uparrow \Lambda }\ E_{\lambda }=E^{\ast },\ $namely%
\begin{equation*}
E^{\ast }:=\left\{ \lim_{\lambda \uparrow \Lambda }\psi (\lambda )\ |\ \psi
(\lambda )\in E\right\} .
\end{equation*}%
$E^{\ast }$ is called the \textbf{natural extension }of $E.$
\end{definition}

Notice that, while the $\Lambda $-limit of a sequence of numbers with
constant value $r\in\mathbb{R}$ is $r$, the $\Lambda$-limit of a constant
sequence of sets with value $E\subseteq\mathbb{R}$ gives a larger set,
namely $E^{\ast }$. In general, the inclusion $E\subseteq E^{\ast }$ is
proper.

This definition, combined with axiom ($\Lambda $-1$)$, entails that 
\begin{equation*}
\mathbb{K}=\mathbb{R}^{\ast }.
\end{equation*}

Given any set $E,$ we can associate to it two sets: its natural extension $%
E^{\ast }$ and the set $E^{\sigma },$ where%
\begin{equation}
E^{\sigma }=\left\{ x^{\ast }\ |\ x\in E\right\} .  \label{sigmaS}
\end{equation}

Clearly $E^{\sigma }$ is a copy of $E;$ however it might be different as a
set since, in general, $x^{\ast }\neq x.$ Moreover $E^{\sigma }\subset
E^{\ast }$ since every element of $E^{\sigma }$ can be regarded as the $%
\Lambda $-limit\ of a constant sequence.

\begin{definition}
\label{limito2}Let 
\begin{equation*}
f_{\lambda }:\ E_{\lambda }\rightarrow \mathbb{R},\ \ \lambda \in \mathfrak{L%
},
\end{equation*}%
be a family of functions. We define a function%
\begin{equation*}
f:\left( \lim_{\lambda \uparrow \Lambda }\ E_{\lambda }\right) \rightarrow 
\mathbb{R}^{\ast }
\end{equation*}%
as follows: for every $\xi \in \left( \lim_{\lambda \uparrow \Lambda }\
E_{\lambda }\right) $ we pose%
\begin{equation*}
f\left( \xi \right) :=\lim_{\lambda \uparrow \Lambda }\ f_{\lambda }\left(
\psi (\lambda )\right) ,
\end{equation*}%
where $\psi (\lambda )$ is a net of numbers such that 
\begin{equation*}
\psi (\lambda )\in E_{\lambda }\ \ \text{and}\ \ \lim_{\lambda \uparrow
\Lambda }\psi (\lambda )=\xi .
\end{equation*}%
A function which is a $\Lambda $-\textit{limit\ is called \textbf{internal}.}
In particular if, $\forall \lambda \in \mathfrak{L,}$ 
\begin{equation*}
f_{\lambda }=f,\ \ \ \ f:\ E\rightarrow \mathbb{R},
\end{equation*}%
we set 
\begin{equation*}
f^{\ast }=\lim_{\lambda \uparrow \Lambda }\ f_{\lambda }.
\end{equation*}%
$f^{\ast }:E^{\ast }\rightarrow \mathbb{R}^{\ast }$ is called the \textbf{%
natural extension }of $f.$
\end{definition}

More in general, the $\Lambda $-limit can be extended to a larger family of
nets. To this aim, let us consider a net%
\begin{equation}
\varphi :\mathfrak{L}\rightarrow {\mathbb{U}}_{n}.  \label{carpa}
\end{equation}%
We will define $\lim\limits_{\lambda \uparrow \Lambda }\varphi (\lambda )$
by induction on $n$. For $n=0,$ $\lim\limits_{\lambda \uparrow \Lambda
}\varphi (\lambda )$ is defined by Th. \ref{nuovo}; so by induction we may
assume that the limit is defined for $n-1$ and we define it for the net (\ref%
{carpa}) as follows:%
\begin{equation}
\lim_{\lambda \uparrow \Lambda }\varphi (\lambda )=\left\{ \lim_{\lambda
\uparrow \Lambda }\psi (\lambda )\ |\ \psi :\mathfrak{L}\rightarrow \mathcal{%
\mathbb{U}}_{n-1}\text{ and}\ \forall \lambda \in \mathfrak{L},\ \psi
(\lambda )\in \varphi (\lambda )\right\} .  \label{limitu}
\end{equation}

\begin{definition}
A mathematical entity (number, set, function or relation) which is the $%
\Lambda $-limit of a net is called \textbf{internal}.
\end{definition}

Let us note that, if $\left( f_{\lambda }\right) $, $\left( E_{\lambda
}\right) $ are, respectively, a net of functions and a net of sets, the $%
\Lambda -$limit of these nets defined by \ref{limitu}) coincides with the $%
\Lambda $-limit given by Definitions \ref{limito} and \ref{limito2}. The
following theorem is a fundamental tool in using the $\Lambda $-limit:

\begin{theorem}
\label{limit}\textbf{(Leibniz Principle)} Let $\mathcal{R}$ be a relation in 
{$\mathbb{U}$}$_{n}$ for some $n\geq 0$ and let $\varphi $,$\psi :\mathfrak{L%
}\rightarrow {\mathbb{U}}_{n}$. If 
\begin{equation*}
\forall \lambda \in \mathfrak{L},\ \varphi (\lambda )\mathcal{R}\psi
(\lambda )
\end{equation*}%
then%
\begin{equation*}
\left( \underset{\lambda \uparrow \Lambda }{\lim }\varphi (\lambda )\right) 
\mathcal{R}^{\ast }\left( \underset{\lambda \uparrow \Lambda }{\lim }\psi
(\lambda )\right) .
\end{equation*}
\end{theorem}

When $\mathcal{R}$ is $\in $ or $\mathcal{=}$ we will not use the symbol $%
\ast $ to denote their extensions, since their meaning is unaltered in
universes constructed over $\mathbb{R}^{\ast }.$ To give an example of how
Leibniz Principle can be used to prove facts about internal entities, let us
prove that if $K\subseteq \mathbb{R}$ is a compact set and $(f_{\lambda })$
is a net of continuous functions then $f=\underset{\lambda \uparrow \Lambda }%
{\lim }f_{\lambda }$ has a maximum on $K^{\ast }$. For every $\lambda $ let $%
\xi _{\lambda }$ be the maximum value attained by $f_{\lambda }$ on $K$, and
let $x_{\lambda }\in K$ be such that $f_{\lambda }(x_{\lambda })=\xi
_{\lambda }.$ For every $\lambda ,$ for every $y_{\lambda }\in K$ we have
that $f_{\lambda }(y_{\lambda })\leq f_{\lambda }(x_{\lambda }).$ By Leibniz
Principle, if we pose 
\begin{equation*}
x=\lim_{\lambda \uparrow \Lambda }x_{\lambda }
\end{equation*}%
we have that%
\begin{equation*}
\forall y\in K\text{ \ }f(y)\leq f(x),
\end{equation*}
so $\xi =\lim_{\lambda \uparrow \Lambda }\xi _{\lambda }$is the maximum of $%
f $ on $K$ and it is attained on $x.$

\subsection{Ultrafunction theory}

Let $\Omega $ be a set in $\mathbb{R}^{N}$ and let $V(\Omega )\ $be a (real
or complex) vector space such that $\mathcal{D}(\overline{\Omega })\subseteq
V(\Omega )\subseteq L^{2}(\Omega )\cap \mathcal{C}(\overline{\Omega }).$

\begin{definition}
Given the function space $V(\Omega )$ we set%
\begin{equation*}
V_{\Lambda}(\Omega ):=\ \underset{\lambda \uparrow \Lambda }{\lim }%
V_{\lambda }(\Omega ),
\end{equation*}%
where%
\begin{equation*}
V_{\lambda }(\Omega )=Span(V(\Omega )\cap \lambda ).
\end{equation*}%
$V_{\Lambda}(\Omega )$ will be called the \textbf{space of ultrafunctions}
generated by $V(\Omega ).$
\end{definition}

Using the above definition, if $V(\Omega )$, $\Omega \subset \mathbb{R}^{N}$%
, is a real function space then we can associate to it three functions
spaces of hyperreal functions, namely $V(\Omega )^{\sigma },$ $%
V_{\Lambda}(\Omega )$ and $V(\Omega )^{\ast }$: 
\begin{equation}
V(\Omega )^{\sigma }=\left\{ f^{\ast }\ |\ f\in V(\Omega )\right\}
\label{sigma1}
\end{equation}%
\begin{equation}
V_{\Lambda}(\Omega )=\left\{ \lim_{\lambda \uparrow \Lambda }\ f_{\lambda }\
|\ f_{\lambda }\in V_{\lambda }(\Omega )\right\}  \label{tilda}
\end{equation}

\begin{equation}
V(\Omega )^{\ast }=\left\{ \lim_{\lambda \uparrow \Lambda }\ f_{\lambda }\
|\ f_{\lambda }\in V(\Omega )\right\}  \label{star}
\end{equation}%
Clearly we have%
\begin{equation*}
V(\Omega )^{\sigma }\subset V_{\Lambda}(\Omega )\subset V(\Omega )^{\ast }.
\end{equation*}

Let us see the relations of the space of ultrafunctions $V_{\Lambda}(\Omega
) $ with the space of ``standard functions'' $V(\Omega )^{\sigma }$(see \ref%
{sigma1}) and the space of internal functions $V(\Omega )^{\ast }$ (see (\ref%
{star})). Given any vector space of functions $V(\Omega )$, the space of
ultrafunction generated by $V(\Omega )$ is a vector space of hyperfinite
dimension that includes $V(\Omega )^{\sigma }$, and the ultrafunctions are $%
\Lambda $-limits of functions in $V_{\lambda }$. Hence the ultrafunctions
are particular internal functions 
\begin{equation*}
u:\left( \mathbb{R}^{\ast }\right) ^{N}\rightarrow {\mathbb{C}^{\ast }.}
\end{equation*}

Since $V_{\Lambda}(\Omega )\subset \left[ L^{2}(\mathbb{R})\right] ^{\ast },$
it can be equipped with the following scalar product%
\begin{equation*}
\left( u,v\right) =\int^{\ast }u(x)\overline{v(x)}\ dx,
\end{equation*}%
where $\int^{\ast }$ is the natural extension of the Lebesgue integral
considered as a functional%
\begin{equation*}
\int :L^{1}(\Omega )\rightarrow {\mathbb{C}}.
\end{equation*}%
Notice that the Euclidean structure of $V_{\Lambda}(\Omega )$ is the $%
\Lambda $-limit of the Euclidean structure of every $V_{\lambda }$ given by
the usual $L^{2}$ scalar product. The norm of an ultrafunction will be given
by 
\begin{equation*}
\left\Vert u\right\Vert =\left( \int^{\ast }|u(x)|^{2}\ dx\right) ^{\frac{1}{%
2}}.
\end{equation*}

\bigskip

\subsection{Morse theory}

Let $\mathfrak{M}$ be a finite dimensional Riemannian manifold and let 
\begin{equation*}
J:\mathfrak{M}\rightarrow \mathbb{R}
\end{equation*}%
be a functional of class $C^{2}.$

A point $u\in \mathfrak{M,}$ is called critical point of $J$ if $dJ(u)=0.$ A
number $c\in \mathbb{R}$ is called critical value of $J$ if there is a
critical point $u\in \mathfrak{M}$ such that $J(u)=c$. A critical point is
called nondegenerate if $H_{J}(u)$ is non singular, namely if%
\begin{equation*}
\left[ \forall \varphi \in T_{u}\mathfrak{M},\ H_{J}(u)\left[ \psi ,\varphi %
\right] =0\right] \Rightarrow \psi =0
\end{equation*}%
If $a,b\in \mathbb{R}$, we set%
\begin{eqnarray*}
J^{b} &=&\left\{ u\in \mathfrak{M\ }|\ \ J(u)\leq b\right\} \\
J_{a}^{b} &=&J^{b}\backslash J^{a}=\left\{ u\in \mathfrak{M\ }|\ \
a<J(u)\leq b\right\} \\
K_{a}^{b} &=&\left\{ u\in J_{a}^{b}\mathfrak{\ }|\ \ dJ(u)=0\right\}
\end{eqnarray*}

The Morse index of a quadratic form $a\left[ \varphi \right] $ is the number
of negative eigenvalues of any matrix representation of $a\left[\varphi %
\right] .$ The Morse index of a critical point $u,$ denoted by $m(u),$ is
the Morse index of the Hessian quadratic form $H_{J}(u)\left[ \varphi \right]
.$ If $u$ is a nondegenerate critical point, we define the \textbf{polynomial%
} \textbf{Morse index }of $u$ as follows%
\begin{equation*}
i_{t}(u)=t^{m(u)}
\end{equation*}%
We have introduced the notion of polynomial Morse index because this notion
allows to define the index of any isolated critical point, even if it is
degenerate; the definition is the following:%
\begin{equation*}
i_{t}(u)=\sum_{k=0}^{N}\dim \left[ H^{k}(J^{c},J^{c}\backslash \left\{
u\right\} )\right] \ t^{k},\ \ \ c=J(u)
\end{equation*}%
where $N$ is the dimension of the manifold $\mathfrak{M,}$ $H^{k}(A,B)$ is
the $k$-th Alexander-Spanier cohomology group of the couple $(A,B)$ with
real coefficients, we denote by $\dim\left[H^{k}(A,B)\right]$ the dimension
of $H^{k}(A,B)$ regarded as real vector space. It is a well known fact of
Morse theory that, if $u$ is a nondegenerate critical point, the two
definitions of $i_{t}(u)$ agree.

We define the Morse polynomial of $J_{a}^{b}$ as follows:%
\begin{equation*}
M_{t}(J_{a}^{b})=\sum_{u\in K_{a}^{b}}i_{t}(u)
\end{equation*}%
Thus $M(t)$ is a polynomial with coefficients in $\mathbb{N}\cup \left\{
+\infty \right\} $. If all the critical points in $K_{a}^{b}$ are not
degenerate, $M(1)$ is the cardinality of $K_{a}^{b}$ namely the number of
the critical points of $J$ in $J_{a}^{b}$. If some critical point is
degenerate, then $M(1)$ is the number of critical points counted with their
multiplicity where the multiplicity of a critical point $u$ is given by $%
i_{1}(u).$

The Betti (or Poincar\'{e}) polynomial of $J_{a}^{b}$ is a topological
invariant defined as follows:%
\begin{equation*}
P_{t}(J_{a}^{b})=\sum_{k=0}^{N}\dim \left[ H^{k}(J^{b},J^{a})\right] \ t^{k}
\end{equation*}
$\dim \left[ H^{k}(J^{b},J^{a})\right] $ is called the $k$-th Betti number
of $J_{a}^{b}.$

In the rest of the paper, we shall use the following important result in
Morse theory.

\bigskip

\begin{theorem}
\label{mara}Let us assume that

\begin{itemize}
\item $\overline{J_{a}^{b}}$ is compact (or more in general $J$ satisfy (PS)
in $\left[ a,b\right] $),

\item $K_{a}^{b}$ is a finite set.
\end{itemize}

Then both $M_{t}(J_{a}^{b})$ and $P_{t}(J_{a}^{b})$ are finite and there
exists a polynomial $Q$ with coefficients in $\mathbb{N}$ such that%
\begin{equation*}
M_{t}(J_{a}^{b})=P_{t}(J_{a}^{b})+(1+t)Q(t)
\end{equation*}
\end{theorem}

\bigskip

\section{Morse theory for ultrafunctions\label{mtu}}

\subsection{Basic results}

Let $V\subset C^{1}(\Omega )$ be a Banach space and let 
\begin{equation*}
J:V\rightarrow \mathbb{R}
\end{equation*}%
be a functional of class $C^{2}$. In the applications, we will assume that $%
J $ has the following structure:%
\begin{equation}
J\left( u\right) =\int F(x,u,\nabla u)\ dx  \label{J}
\end{equation}

As we emphasized in the introduction the main difficult for the development
of Morse Theory in Banach spaces is to define the right concept of
nondegeneracy and of Morse index for a critical point.

We will be interested in Morse theory for the functional 
\begin{equation*}
J_{\Lambda }:V_{\Lambda }\rightarrow \mathbb{R}^{\ast }
\end{equation*}%
where $V_{\Lambda }$ is a space of ultrafunctions and $J_{\Lambda }$ is the
restriction of $J^{\ast }$ to $V_{\Lambda }.$For example, a suitable space
for the functional (\ref{J}) is $V_{\Lambda }(\Omega ):=[C^{2}(\Omega )\cap
C_{0}^{1}(\overline{\Omega })]_{\Lambda }.$

Now let us describe the main objects of Morse theory in the ultrafunctions
framework.

\begin{definition}
An ultrafunction $u\in V_{\Lambda}$ is called a critical point of $%
J_{\Lambda}:V_{\Lambda}\rightarrow \mathbb{R}^{\ast }$ if%
\begin{equation*}
\forall \varphi \in V_{\Lambda},\ d J_{\Lambda}(u)\left[ \varphi \right] =0
\end{equation*}%
where $dJ$ is the differential of $J.$
\end{definition}

In particular, if $J$ is the functional (\ref{J}), we have that $u\in
V_{\Lambda }=[C^{2}(\Omega )\cap C_{0}^{1}(\overline{\Omega })]_{\Lambda }$
is a critical point if%
\begin{equation*}
\forall \varphi \in V_{\Lambda }(\Omega ),\ \int \left[ \frac{\partial F}{%
\partial \left( \nabla u\right) }\cdot \nabla \varphi +\frac{\partial F}{%
\partial u}\varphi \right] \ dx=0
\end{equation*}%
Here $\frac{\partial F}{\partial \left( \nabla u\right) }$ denotes the
vector $\left( \frac{\partial F}{\partial u_{x_{1}}},....,\frac{\partial F}{%
\partial u_{x_{N}}}\right) .$

The Hessian quadratic form $H_{J^{\ast }}(u)$ of $J^{\ast }$ is defined on $%
V^{\ast }\times V^{\ast };$ we will denote by $H_{J_{\Lambda }}(u)$ its
restriction to $V_{\Lambda }\times V_{\Lambda }.$ A critical point of $%
J_{\Lambda }$ is called nondegenerate if%
\begin{equation*}
\forall \varphi \in V_{\Lambda },\ H_{J_{\Lambda }}(u)\left[ \psi ,\varphi %
\right] =0\Rightarrow \psi =0
\end{equation*}%
Since $H_{J_{\Lambda }}(u)$ is a quadratic form defined on a hyperfinite
space $V_{\Lambda },$ its Morse index is well defined and hence also the
Morse index $m_{\Lambda }(u)$ of $u$ is well defined.

Given two hyperreal numbers $a<b,$ we set%
\begin{align*}
J^{b}_{\Lambda} &=\left\{ u\in V_{\Lambda}\mathfrak{\ }|\ \
J_{\Lambda}(u)\leq b\right\} \\
[J_{a}^{b}]_{\Lambda}&=J_{\Lambda}^{b}\backslash J^{a}_{\Lambda}=\left\{
u\in V_{\Lambda}\mathfrak{\ }|\ \ a<J_{\Lambda}(u)\leq b\right\} \\
[K_{a}^{b}]_{\Lambda} &=\left\{ u\in J_{a}^{b}\mathfrak{\ }|\ \
dJ_{\Lambda}(u)=0\right\}
\end{align*}

Next we must define the Morse index, the Morse polynomial and the Betti
polynomial in the frame of ultrafunctions. We could define them
intrinsically as we have done for the above notions. However it seems easier
to define them by mean of a $\Lambda $-limit.

We set%
\begin{equation*}
M_{t}([J_{a}^{b}]_{\Lambda})=\lim_{\lambda \uparrow \Lambda }\
M_{t}(J_{a_{\lambda }}^{b_{\lambda }}\cap V_{\lambda })
\end{equation*}%
where $a_{\lambda }\ $and $b_{\lambda }$ are two real nets such that%
\begin{equation}
\lim_{\lambda \uparrow \Lambda }\ a_{\lambda }=a;\ \ \lim_{\lambda \uparrow
\Lambda }\ b_{\lambda }=b.  \label{lin}
\end{equation}

Analogously, we define the "generalized" Betti polynomial as follows:%
\begin{equation*}
P_{t}([J_{a}^{b}]_{\Lambda})=\lim_{\lambda \uparrow \Lambda }\
P_{t}(J_{a_{\lambda }}^{b_{\lambda }}\cap V_{\lambda }).
\end{equation*}

Now it is possible to state an abstract theorem for Morse theory in the
framework of ultrafunctions:

\begin{theorem}
\label{cecilia}Let 
\begin{equation*}
J:V\rightarrow \mathbb{R}
\end{equation*}%
be a $C^{2}$-functional and 
\begin{equation*}
J_{\Lambda}:V_{\Lambda}\rightarrow \mathbb{R}^{\ast }
\end{equation*}%
be the restriction of $J^{\ast }$ to $V_{\Lambda}.$ Let $a,b\in \mathbb{R}%
^{\ast }$ satisfy (\ref{lin}) and assume that

\begin{itemize}
\item for almost every $\lambda \in \mathfrak{L},$ $\overline{J_{a_{\lambda
}}^{b_{\lambda }}}$ is compact (or more in general $J$ satisfy (PS) in $%
\left[ a_{\lambda },b_{\lambda }\right] $),

\item for almost every $\lambda \in \mathfrak{L},\mathfrak{\ }K_{a_{\lambda
}}^{b_{\lambda }}\ $is finite .
\end{itemize}

Then $M_{t}([J_{a}^{b}]_{\Lambda}),P_{t}([J_{a}^{b}]_{\Lambda})\in \mathfrak{%
pol}(\mathbb{N}\mathbf{)}^{\ast }$ where%
\begin{equation*}
\mathfrak{pol}(\mathbb{N}\mathbf{)=}\left\{ \text{polynomials\ with
coefficients in }\mathbb{N}\right\}
\end{equation*}%
and there exists a polynomial $Q\in \mathfrak{pol}(\mathbb{N}\mathbf{)}%
^{\ast }$ such that%
\begin{equation*}
M_{t}([J_{a}^{b}]_{\Lambda})=P_{t}([J_{a}^{b}]_{\Lambda})+(1+t)Q(t).
\end{equation*}
\end{theorem}

\textbf{Proof - }For almost every $\lambda \in \mathfrak{L},$ $\overline{%
J_{a_{\lambda }}^{b_{\lambda }}}$ is compact and $K_{a_{\lambda
}}^{b_{\lambda }}\ $is finite; then by Th. \ref{mara}, $M_{t}(J_{a_{\lambda
}}^{b_{\lambda }})$ and $P_{t}(J_{a_{\lambda }}^{b_{\lambda }})\in \mathfrak{%
pol}(\mathbb{N}\mathbf{)}$ and there exists a polynomial $Q_{\lambda }\in 
\mathfrak{pol}(\mathbb{N}\mathbf{)}$ such that%
\begin{equation*}
M_{t}(J_{a_{\lambda }}^{b_{\lambda }})=P_{t}(J_{a_{\lambda }}^{b_{\lambda
}})+(1+t)Q_{\lambda }(t)
\end{equation*}%
The conclusion follows taking the $\Lambda $-limit.

$\square $

\bigskip

\subsection{Ultrafunctions versus Sobolev spaces}

\bigskip

Usually, the critical points of functional of type (\ref{J}) are studied in
the Sobolev space $W_{0}^{1,p}(\Omega )$ provided that the functional $J$
can be extended to $W_{0}^{1,p}(\Omega )$ as a $C^{1}$ functional. In this
section, we will investigate some relation between the ultrafunction and the
Sobolev space approach.

So we will assume that $J$ can be extended to a $C^{1}$-functional in a
Banach space $W\subset L^{1}(\Omega )$ (with some abuse of notation we will
denote this extension by the same letter $J$):%
\begin{equation*}
J:W\rightarrow \mathbb{R}.
\end{equation*}%
So, we have that%
\begin{equation*}
V^{\sigma }\subset W^{\sigma }\subset V_{\Lambda}
\end{equation*}

In the following, to simplify the notation, we will identify $V^{\sigma }$
and $V$ as well as $W^{\sigma }$ and $W.$

\bigskip

The next theorems will establish some relations between the critical points
of $J_{\Lambda}$ in $V_{\Lambda}$ and the critical points of $J$ in $W.$

The first result in this direction is (almost) trivial:

\begin{theorem}
Under the same framework and the same assumptions of Th. \ref{cecilia} every
critical point of $J$ in $W$ is a critical point of $J_{\Lambda }$ in $%
V_{\Lambda }$
\end{theorem}

\textbf{Proof: }Let $u\in W$ be a critical point of $J$; we will use the
fact that $V(\Omega )^{\sigma }\subset V_{\Lambda }(\Omega )$ to prove the
thesis.

Let $u_{\lambda }$ be the constant net $u_{\lambda }=u$; then 
\begin{equation*}
\lim_{\lambda \uparrow \Lambda }u_{\lambda }=u^{\ast }\in V(\Omega )^{\sigma
}\subset V_{\Lambda }(\Omega ),
\end{equation*}%
and let $J_{\lambda }$ be the constant net $J_{\lambda }=J$; then 
\begin{equation*}
dJ_{\lambda }(u_{\lambda })[\phi _{\lambda }]=0
\end{equation*}%
for every $\phi _{\lambda }\in V_{\lambda }(\Omega )$; therefore, taking the 
$\Lambda $-limit of a constant net we have the thesis.

$\square $

The above theorem cannot be inverted in the sense that it is false that
every critical point of $J_{\Lambda}$ is a critical point of $J\ $in$\ W.$
However, there are conditions which insure the existence of critical point
of $J$ in $W$. More precisely the next theorem states that, under suitable
condition, ``infinitely close'' to any critical point of $J_{\Lambda}$ there
is is a critical point of $J$

This theorem exploit a compactness condition which is a variant of the usual
Palais-Smale condition (PS). We recall the Palais-Smale condition is a basic
tool for Morse theory in infinite dimensional manifolds (see e.g. \cite%
{Chang93}). Here it is used only to relate some critical point of $%
J_{\Lambda }$ with the critical points of $J.$

\begin{definition}
\label{PSU}\textbf{Palais-Smale condition for ultrafunctions (PSU)} We say
that the functional%
\begin{equation*}
J:W\rightarrow \mathbb{R}
\end{equation*}%
satisfies (PSU) in the interval $\left[ a,b\right] \subset \mathbb{R}$ if
for every net $\left\{ u_{\lambda }\right\} _{\lambda \in \mathfrak{L}}$
such that that

\begin{itemize}
\item (A) $\forall \lambda \in \mathfrak{L,\ }J(u_{\lambda })\in \left[ a,b%
\right] $

\item (B) $\forall \lambda \in \mathfrak{L},\ \forall v\in V_{\lambda },\
dJ(u_{\lambda })\left[ v\right] =0$
\end{itemize}

\noindent there is a converging subnet $\left\{ u_{\lambda }\right\}
_{\lambda \in \mathfrak{D}}$ ($\mathfrak{D}\subset \mathfrak{L}$) in the
topology of $W$, such that 
\begin{equation*}
\lim_{\lambda \rightarrow \Lambda }\ u_{\lambda }\in W.
\end{equation*}
\end{definition}

\bigskip

\begin{remark}
Notice that, by prop. \ref{nino}, the sequence $\left\{ u_{\lambda }\right\}
_{\lambda \in \mathfrak{L}}$ itself is converging.
\end{remark}

\begin{theorem}
\label{A}Let us assume that $W$ is a Banach space and that $V\subset
W\subset V_{\Lambda }.$ Let 
\begin{equation*}
J:W\rightarrow \mathbb{R}
\end{equation*}%
be a $C^{1}$-functional which satisfies (PSU) in the interval $\left[ a,b%
\right] .\ $Then, if $\bar{u}$ is a critical point of 
\begin{equation*}
J_{\Lambda }:V_{\Lambda }\rightarrow \mathbb{R}^{\ast }
\end{equation*}%
with $J_{\Lambda }\left( \bar{u}\right) \in \left[ a,b\right] ^{\ast },$
there exists $w\in K_{a}^{b}$ such that%
\begin{equation*}
\left\Vert \bar{u}-w^{\ast }\right\Vert _{W^{\ast }}\sim 0.
\end{equation*}
\end{theorem}

\begin{remark}
Notice that in the above theorem, it is possible that that $\bar{u}=w^{\ast
}.$ Obviously, this fact always occur if $W$ is a Hilbert space and all the
critical values of $J$ in $\left[ a,b\right] $ are not degenerate.
\end{remark}

\textbf{Proof of Th. \ref{A}. }Let%
\begin{equation*}
\bar{u}=\lim_{\lambda \uparrow \Lambda }\ u_{\lambda };
\end{equation*}
Then, since (PSU) holds, there is a function $w\in W$ and a subnet of $%
u_{\lambda }$ such that 
\begin{equation*}
\left\Vert u_{\lambda }-w\right\Vert _{W}\rightarrow 0.
\end{equation*}%
By Proposition \ref{nino}, $\left\Vert u_{\lambda }-w\right\Vert $ is a
converging net, and hence, for every $\varepsilon >0,$ exists $Q\in \mathcal{%
U}$ such that $\forall \lambda \in Q,\ $ 
\begin{equation*}
\left\Vert u_{\lambda }-w\right\Vert _{W}\leq \varepsilon .
\end{equation*}%
If you take the $\Lambda $-limit of the above inequality, you get that%
\begin{equation*}
\left\Vert \bar{u}-w^{\ast }\right\Vert _{W^{\ast }}\leq \varepsilon .
\end{equation*}%
By the arbitrariety of $\varepsilon $, we conclude that 
\begin{equation*}
\left\Vert \bar{u}-w^{\ast }\right\Vert _{W^{\ast }}\sim 0
\end{equation*}

$\square $


\bigskip

\begin{thebibliography}{99}
%\bibitem{AmRa1973} A. Ambrosetti and P.H. Rabinowitz, \emph{Dual Variational
%Methods in Critical Point Theory and Applications,} J. Funct. Anal., Vol.
%14, (1973), pp. 349-381.

%\bibitem{bahri} Bahri A., \textsl{Critical points at infinity in the
%variational calculus, }in:\textsl{\ } Partial differential equations (Rio de
%Janeiro, 1986), Lecture Notes in Math., \textbf{1324}, Springer, Berlin,
%(1988), p. 1-29.

%\bibitem{benci95} Benci V., \textsl{A construction of a nonstandard universe}%
%, in: Advances of Dynamical Systems and Quantum Physics (S. Albeverio et
%al., eds.), World Scientific, Singapore, (1995), p. 11-21.

%\bibitem{benci99} Benci V., \textsl{An algebraic approach to nonstandard
%analysis, }in: Calculus of Variations and Partial differential equations,%
%\textsl{\ }p. 285-307, (G.Buttazzo, et al., eds.), Springer, Berlin, (1999),
%p.285-326.

\bibitem{ultra} Benci V., \textsl{Ultrafunctions and generalized solutions,}
in: Adv. Nonlinear Stud. 13, (2013), 461--486, arXiv:1206.2257.

%\bibitem{BGG} Benci V., Galatolo S., Ghimenti M., \textsl{An elementary
%approach to Stochastic Differential Equations using the infinitesimals}, in
%Contemporary Mathematics, \textbf{530}, Ultrafilters across Mathematics,
%American Mathematical Society, (2010), p. 1-22.

%\bibitem{BDN2003} Benci V., Di Nasso M., \textsl{Alpha-theory: an elementary
%axiomatic for nonstandard analysis}, Expositiones Mathematicae \textbf{21}
%(2003) p.~355--386.

%\bibitem{BHW} Benci V., Horsten H., Wenmackers S., \textsl{Non-Archimedean
%probability}, submitted, URL: http://arxiv.org/abs/1106.1524.

\bibitem{belu2013} Benci V., Luperi Baglini L., \textsl{Basic Properties of
ultrafunctions,} to appear in the WNDE2012 Conference Proceedings,
arXiv:1302.7156.

\bibitem{milano} Benci V., Luperi Baglini L., \textsl{Ultrafunctions and
applications}, to appear on DCDS-S (Vol. 7, No. 4) August 2014,
arXiv:1405.4152.

\bibitem{beyond} Benci V., Luperi Baglini L., \textsl{Generalized functions
beyond distributions}, to appear on AJOM (2014), arXiv:1401.5270.

\bibitem{topology} Benci V., Luperi Baglini L., \textsl{A topological
approach to Non Archimedean Mathematics}, to appear, arXiv:1412.2223

%\bibitem{BN} Brezis H., Nirenberg L., \textsl{Positive solutions of
%nonlinear elliptic equations involving criticalSobolev exponents,}
%Comm. Pure Appl. Math., 36 (1983), p. 437--477

\bibitem{CCV03} J. Carmona, S. Cingolani, G. Vannella, \textsl{Estimates of
the critical groups for solutions of quasilinear elliptic systems},
Electron. J. Differential Equations (2003), 1--13.

\bibitem{Cha} Chabrowski J., \textsl{Variational methods for potential
operator equations, with applications to nonlinear elliptic equations.}
Walter de Gruyter \& Co., Berlin, 1997.

\bibitem{Chang83} K. Chang \textsl{Morse theory on Banach space and its
applications to partial differential equations} Chin. Ann. of Math., 4B
(1983), pp. 381-399

\bibitem{Chang93} K. Chang \textsl{Infinite Dimensional Morse Theory and
Multiple Solution Problems} Birkh\~{A}\textcurrency user, Boston (1993)

\bibitem{Chang98} K. Chang \textsl{Morse theory in nonlinear analysis} in A.
Ambrosetti, K.C. Chang, I. Ekeland (Eds.), Nonlinear Functional Analysis and
Applications to Differential Equations, Word Scientific, Singapore (1998)

\bibitem{CinCarMartVan} S. Cingolani, J. Carmona P. J. Mart\TEXTsymbol{%
\backslash}'\{i\}nez-Aparicio, G. Vannella, \textsl{Regularity and Morse
index of the solutions to critical Quasilinear Elliptic Systems},
Communications in Partial Differential Equations, vol. 38 (2013)

\bibitem{CD05} S. Cingolani, M. Degiovanni, \textsl{Nontrivial solutions for
p-Laplace equations with right hand side having p-linear growth at infinity}%
, Comm. Partial Differential Equations, vol. 30 (2005), 1191--1203.

\bibitem{CinDe09} S. Cingolani, M. Degiovanni, \textsl{On the Poincar%
\TEXTsymbol{\backslash}'\{e\}-Hopf Theorem for Functionals defined on Banach
Spaces}, Advances Nonlinear Studies, vol. 9 (2009), 679--699.

\bibitem{CinDeSciCGpLap} S. Cingolani, M. Degiovanni, B. Sciunzi, \textsl{%
Critical groups estimates for \$p\$-Laplace equations via Uniform Sobolev
Inequalities}, submitted for publication.

\bibitem{CLV05} S. Cingolani, M. Lazzo, G. Vannella, \textsl{Multiplicity
results for a quasilinear elliptic system via Morse theory}, Communications
in Contemporary Mathematics, vol.7 (2005), 227--249.

\bibitem{CV03} S. Cingolani, G. Vannella, \textsl{Critical groups
computations on a class of Sobolev Banach spaces via Morse index}, Ann.
Inst. H. Poincar\TEXTsymbol{\backslash}`\{e\} Anal. Non Lineaire, vol. 2
(2003), 271-292.

\bibitem{CV06} S. Cingolani, G. Vannella, \textsl{Morse index and critical
groups for p-Laplace equations with critical exponents}, Mediterranean
Journal of Mathematics , vol.3 (2006), 347--592.

\bibitem{CV07} S. Cingolani, G. Vannella, \textsl{Marino-Prodi perturbation
type results and Morse indices of minimax critical points for a class of
functionals in Banach spaces}, Annali di Matematica Pura e Applicata, vol.
186 (2007), 157--185.

\bibitem{CV09} S. Cingolani, G. Vannella, \textsl{Multiple positive
solutions for a critical quasilinear equation via Morse theory}, Ann. Inst.
H. Poincar\TEXTsymbol{\backslash}'\{e\} Anal. Non Lineaire, vol. 26 (2009),
397--413.

\bibitem{CV09b} S. Cingolani, G. Vannella, \textsl{On the multiplicity of
positive solutions for p-Laplace equations via Morse theory}, Journal of
Differential Equations, vol. 247 (2009), 3011--3027.

\bibitem{CinVanVis11} S. Cingolani, G. Vannella, D. Visetti, \textsl{Morse
index estimates for quasilinear equations on Riemannian manifolds}, Advances
in Differential Equations, vol. 16 (2011), 1001-1020.

\bibitem{CinVanVisMul} S. Cingolani, G.Vannella, D. Visetti, \textsl{%
Multiplicity of positive solutions for a quasilinear equation on a
Riemannian manifold}, submitted for publication.

\bibitem{el06} Ehrlich Ph., \textsl{The Rise of non-Archimedean Mathematics
and the Roots of a Misconception I: The Emergence of non-Archimedean Systems
of Magnitudes,} Arch. Hist. Exact Sci. 60 (2006) 1--121, Identifier (DOI)
10.1007/s00407-005-0102-4.

\bibitem{GroMe} Gromoll D., Meyer W., \textsl{On differentiable functions
with isolated critical points }, Topology Volume 8, Issue 4, 1969, Pages
361-369

\bibitem{keisler76} Keisler H.J., \textsl{Foundations of Infinitesimal
Calculus}, Prindle, Weber \& Schmidt, Boston 1976. [This book is now freely
downloadable at: http://www.math.wisc.edu/\symbol{126}%
keisler/foundations.html]

\bibitem{Lan} S. Lancelotti \textsl{Morse index estimates for continuous
functionals associated with quasilinear elliptic equations} Adv.
Differential Equations, 7 (2002), pp. 99-128

\bibitem{MerPalm} F. Mercuri, G. Palmieri \textsl{Problems in extending
Morse theory to Banach spaces} Boll. UMI, 12 (1975), pp. 397-401

\bibitem{Uhl72} K. Uhlenbeck \textsl{Morse theory on Banach manifolds} J.
Funct. Anal., 10 (1972), pp. 430-445
\end{thebibliography}
\end{document}